\documentclass[a4paper,12pt]{article}
\usepackage{fullpage}
\usepackage[english,francais]{babel}
\usepackage{epsf}
\usepackage[dvips]{epsfig}
\usepackage{amssymb}
\usepackage{amsmath}
\usepackage[all]{xy}
\usepackage{appendix}
\usepackage{fancyhdr,graphicx,color}

\usepackage{ucs}
\usepackage[utf8x]{inputenc}\usepackage{wasysym}
\usepackage{aeguill}

\newtheorem{theorem}{Theorem}
\newtheorem{prop-f}[theoreme]{Proposition}
\newtheorem{prop}[theorem]{Proposition}

\newtheorem{corollary}[theorem]{Corollary}

\newtheorem{lemma}[theorem]{Lemma}

\newcommand{\finpreuve}{\hspace{\stretch{1}}{$\square$}}

\newcommand{\card}{{\rm card}}

\newcommand{\R}{\mathbb{R}}

\newcommand{\Z}{\mathbb{Z}}

\renewcommand{\epsilon}{\varepsilon}
\def\btab{\begin{eqnarray*}}
\def\etab{\end{eqnarray*}}
\def\beq{\begin{equation}}
\def\eeq{\end{equation}}

\newcounter{numeroexo}

\def\as{\overset{\Sigma}{\longleftrightarrow}}

\title{\LARGE Positivity of the time constant in a continuous model of first passage percolation}
\author{Jean-Baptiste Gou\'er\'e\footnote{LMPT, Universit\'e Fran\c cois Rabelais, Parc de Grandmont, 37200 Tours, France, {\it jean-baptiste.gouere@lmpt.univ-tours.fr}} and Marie Th\'eret\footnote{LPMA, Universit\'e Paris Diderot, 5 rue Thomas Mann, 75205 Paris Cedex 13, France, {\it marie.theret@univ-paris-diderot.fr}}}
\date{}

\begin{document}

\selectlanguage{\english}

\maketitle

\thispagestyle{empty}

\noindent
{\bf Abstract:} {We consider a non trivial Boolean model $\Sigma$ on $\R^d$ for $d\geq 2$. For every $x,y \in \R^d$ we define $T(x,y)$ as the minimum time needed to travel from $x$ to $y$ by a traveler that walks at speed $1$ outside $\Sigma$ and at infinite speed inside $\Sigma$. By a standard application of Kingman sub-additive theorem, one easily shows that $T(0,x)$ behaves like $\mu \|x\|$ when $\|x\|$ goes to infinity, where $\mu$ is a constant named the time constant in classical first passage percolation. In this paper we investigate the positivity of $\mu$. More precisely, under an almost optimal moment assumption on the radii of the balls of the Boolean model, we prove that $\mu>0$ if and only if the intensity $\lambda$ of the Boolean model satisfies $\lambda <   \widehat{\lambda}_c$, where $ \widehat{\lambda}_c$ is one of the classical critical parameters defined in continuum percolation.}
\\





\section{Introduction and main results}

\subsection{Boolean model}

The Boolean model is defined as follows.
At each point of a homogeneous Poisson point process on the Euclidean space $\R^d$, we center a ball of random radius.
We assume that the radii of the balls are independent, identically distributed and independent of the point process.
The Boolean model is the union of the balls.
There are three parameters: 
\begin{itemize}
\item An integer $d \ge 2$. This is the dimension of the ambient space $\R^d$.
\item A real number $\lambda>0$.  The intensity measure of the Poisson point process of centers is $\lambda |\cdot|$ where $|\cdot|$ denotes the Lebesgue measure on $\R^d$.
\item A probability measure $\nu$ on $(0,+\infty)$. This is the common distribution of the radii.
\end{itemize}
We will denote the Boolean model by $\Sigma(\lambda,\nu,d)$ or $\Sigma$.

More precisely, the Boolean model is defined as follows. 
Let $\xi$ be a Poisson point process on $\R^d \times (0,+\infty)$ with intensity measure $\lambda|\cdot| \otimes \nu$.
Set
$$
\Sigma(\lambda,\nu,d)=\bigcup_{(c,r) \in \xi} B(c,r)
$$
where $B(c,r)$ denotes the open Euclidean ball of $\R^d$ with center $c$ and radius $r$.
We refer to the book by Meester and Roy \cite{Meester-Roy-livre} for background on the Boolean model, and to the book by Schneider and Weil \cite{Schneider-Weil} and the book by Last and Penrose \cite{LastPenrose} for background on Poisson processes.
We also denote by $S(c,r)$ the Euclidean sphere of $\R^d$ with center $c$ and radius $r$.
We write $S(r)$ when $c=0$.

\bigskip

In this paper, we will always assume  
\begin{equation}\label{e:momentd}
 \int_{(0,+\infty)} r^d \nu(dr) < \infty.
\end{equation}
When \eqref{e:momentd} does not hold, all the models we consider are trivial.
This is due to the fact that, if \eqref{e:momentd} does not hold, then for any $\lambda>0$, with probability one, $\Sigma = \R^d$.
This is Proposition 3.1 in \cite{Meester-Roy-livre}.

Let us state a simple consequence of \eqref{e:momentd}. 
With probability one, the number of random balls which touch a given bounded subset of $\R^d$ is finite
\footnote{Let $n \ge 1$. The number $N_n$ of random balls which touch $B(0,n)$ is a Poisson random variable with parameter 
$$
\lambda \int_{(0,+\infty)} v_d (r+n)^d \mu(dr)
$$
where $v_d$ is the volume of the unit ball of $\R^d$.
Therefore, with probability one, all the $N_n$ are finite.}.

\bigskip

Let $\chi$ denote the set of centers, that is the projection of $\xi$ on $\R^d$.
This is a Poisson point process of intensity measure $\lambda |\cdot|$.
For each $c \in \chi$, we denote by $r(c)$ the unique\footnote{
Consider the projection from $\R^d \times (0,+\infty) \to \R^d$.
With probability one, the restriction to $\xi$ of this projection is one-to-one.} real $r$ such that $(c,r)$ belongs to $\xi$.
When $c \in \R^d \setminus \chi$, we set $r(c)=0$.

\subsection{Paths}

In this paper we only consider polygonal paths. A path is a finite sequence of {\em distinct} points of $\R^d$ - if the points are not distinct, we simply name it a sequence.
The length of a path $\pi=(x_0,\dots,x_k)$ is 
$$
\ell(\pi)=\sum_{i=1}^k \|x_i-x_{i-1}\|
$$
where $\|\cdot\|$ denotes the usual Euclidean norm on $\R^d$.
In some cases, we will also see $\pi$ as a curve $[0,\ell(\pi)] \to \R^d$ parametrized by arc length.
A path from $A \subset \R^d$ to $B \subset \R^d$ is a path such that $\pi(0) \in A$ and $\pi(\ell(\pi)) \in B$.
A path is in $C \subset \R^d$ if $\pi([0,\ell(\pi)]) \subset C$. Notice that if $\pi=(x_0,\dots,x_k)$ is a path, then $\pi([0,\ell(\pi)])  $ is the finite union of the closed segments $[x_{i-1}, x_i]$ for $i\in \{1,\dots , k\}$. All these definitions can be extended to sequences in a natural way. 

\subsection{Percolation in the Boolean model}

\paragraph{Two critical thresholds.} If $A$ and $B$ are two subsets of $\R^d$, we set
$$
\{A \overset{\Sigma}{\longleftrightarrow} B\} = \{\mbox{There exists a path in }\Sigma\mbox{ from }A\mbox{ to }B\}
$$
and
$$
\{0 \as \infty\} = \{\mbox{The connected component of }\Sigma \mbox{ that contains the origin is unbounded}\}.
$$
We define two critical thresholds by
$$
\lambda_c=\lambda_c(\nu,d) = \sup \{\lambda > 0 : P(0 \as \infty) = 0\} \in [0,+\infty]
$$
and
$$
\widehat{\lambda}_c = \widehat{\lambda}_c(\nu,d) = \sup \{\lambda > 0 :  \lim_{r\rightarrow \infty } P(S(r) \as S(2r)) = 0\} \in [0,+\infty]
$$
where $S(u)$, $u >0$, denotes the Euclidean sphere of radius $u$ centered at the origin.

\bigskip

\paragraph{Non triviality of the thresholds.} Recall that we assume \eqref{e:momentd}. 
The thresholds are non trivial. More precisely,
$$
0 < \widehat{\lambda}_c \le \lambda_c < \infty.
$$
The inequality $\lambda_c < \infty$ is proven for a more general model by Hall in \cite{Hall-continuum-percolation} (see Theorem 3).
In our setting, this can be proven in a simple way by coupling the Boolean percolation model with a Bernoulli percolation model on $\Z^d$.
This is explained in the remark below the proof of Theorem 3.3 in the book by Meester and Roy \cite{Meester-Roy-livre}.
The inequality $\widehat{\lambda}_c \le \lambda_c$ is a consequence of the following simple fact:
$$
P(0 \as \infty) = \lim_{r \to \infty} P(0 \as S(2r)) \le \limsup_{r \to \infty} P(S(r) \as S(2r)) .
$$
The proof of the inequality $0 < \widehat{\lambda}_c$ is implicit  in \cite{G-perco-boolean-model} where one of the main aims is to prove the positivity of $\lambda_c$.
We refer to Appendix \ref{s:explicitons} for more details. 

\paragraph{The set $\{\lambda>0 :  \lim_{r \rightarrow \infty }P(S(r) \as S(2r)) = 0\}$ is open.} This result is implicit in \cite{G-perco-boolean-model}.
We refer to Appendix \ref{s:explicitons} for more details. 

\paragraph{Phase transition.} In particular, $\lambda_c$ is non trivial.
Therefore, there exists a subcritical phase and a supercritical phase for percolation.
\begin{itemize}
 \item If $\lambda < \lambda_c$, then with probability one there is no unbounded component in $\Sigma$.
 \item If $\lambda > \lambda_c$, then with probability one there exists at least one (and actually a unique) unbounded component in $\Sigma$.
\end{itemize}
We refer to \cite{Meester-Roy-livre} and to the book by Penrose \cite{Penrose} for background on percolation in the Boolean model.

\paragraph{Sharp threshold.} The critical parameter $\lambda_c$ is probably the more intuitive to define, however in what follows the relevant critical parameter to consider is $\widehat \lambda_c$. For this reason we present here known results concerning the link between $\lambda_c$ and $\widehat \lambda_c$.

If the radii are bounded, then $\lambda_c=\widehat\lambda_c$.
This is a sharp threshold property.
The sharpness of the transition in the discrete setting was proved independently 
by Menshikov \cite{Menshikov-coincidence} and by Aizenman-Barsky \cite{Aizenman-Barsky}.
The first proof of the equality $\lambda_c=\widehat\lambda_c$ relied on the analogous result in the discrete setting.
We refer to \cite{Meester-Roy-livre} for the proof (see Theorem 3.5) and references. Ziesche gives in \cite{Ziesche} a short proof of the equality $\lambda_c = \widehat \lambda_c$ for bounded radii. It relies on a new and short proof of the analogous result in the discrete setting by Duminil-Copin and Tassion \cite{DCT-sharp-Bernoulli, DCT-sharp-Ising}.

In dimension $2$, the sharpness of the transition is one of the results proven recently by Ahlberg, Tassion and Teixera in \cite{ATT}, using a strategy which is specific to the dimension $2$.

\subsection{First-passage percolation in the Boolean model}

In \cite{GM-deijfen}, R\'egine Marchand and the first author studied a model introduced by Deijfen in \cite{Deijfen-modele}.
The model we introduce in this paper appears implicitly in \cite{GM-deijfen} as an intermediate model.
We refer to \cite{GM-deijfen} for the definition of Deijfen's model and its links with the model defined here.

A traveler walks on $\R^d$.
Inside the Boolean model $\Sigma$ he walks at infinite speed.
Outside the Boolean model $\Sigma$ he walks at speed $1$.
He travels from $x \in \R^d$ to $y \in \R^d$ as fast as he can.
We denote by $T(x,y)$ the time needed to perform this travel.
For example if $x$ and $y$ belong to the same connected component of $\Sigma$, then $T(x,y)=0$.

Here is a more formal definition.
For any $a$ and $b$ in $\R^d$, we define $\tau(a,b)$ as the one-dimensional Hausdorff measure of $[a,b] \cap \Sigma^c$.
With each path $\pi=(x_0,...,x_n)$ is associated a time as follows:
$$
\tau(\pi)=\sum_{i=1}^{n} \tau(x_{i-1},x_i).
$$
If $x$ and $y$ are two points of $\R^d$, then $T(x,y)$ is defined by:
$$
T(x,y)=\inf\{\tau(\pi) : \pi \in {\mathcal C}(x,y)\},
$$
where ${\mathcal C}(x,y)$ is the set of paths from $x$ to $y$.

A standard application of Kingman sub-additive theorem yields the following result.

\begin{theorem}\label{t:forme}
There exists a constant $\mu=\mu(\lambda,\nu,d) \in [0,1]$ such that:
$$
\lim_{\|x\|\to\infty}\frac{T(0,x)}{\|x\|} = \mu \mbox{ with probability }1\mbox{ and in }L^1.
$$
\end{theorem}
We emphasize the fact that the convergence stated in Theorem \ref{t:forme} is uniform in all directions. Note that by the isotropy of the model the asymptotic behavior of $T(0,x)$ does not depend on the direction. Theorem \ref{t:forme} can be deduced from Theorem 1 in Ziesche \cite{Ziesche2}, but we give a proof of this result in Appendix \ref{s:preuve-forme} for self-containedness. For any $A, B \subset \R^d$ we write
$$
T(A,B) = \inf_{a \in A, b \in B} T(a,b).
$$
For any $r>0$, we use the shorthand notation 
$$
T(r)=T(\{0\},S(r)).
$$
By Theorem \ref{t:forme} we get
\begin{equation}\label{e:forme}
\lim_{r\rightarrow \infty} \frac{T(r)}r = \mu \mbox{ a.s. and in  }L^1.
\end{equation}

\subsection{Link between percolation and first passage percolation; main result}

Consider the following condition:
\begin{equation}\label{e:greedy}
 \int_{(0,+\infty)} \nu([r,+\infty))^{1/d} dr < \infty .
\end{equation}
It appears in the paper by Martin \cite{Martin-greedy} about greedy lattice paths and animals. 
We refer to \cite{Martin-greedy} for a discussion about Condition \eqref{e:greedy}.
For example, for any $\epsilon>0$,
$$
\int_{(0,+\infty)} r^d \ln_+(r)^{d-1+\epsilon}\nu( dr) < \infty 
\Rightarrow \int_{(0,+\infty)} \nu([r,+\infty))^{1/d} dr < \infty
\Rightarrow \int_{(0,+\infty)} r^d \nu (dr) < \infty.
$$

Here is the main result of this paper.

\begin{theorem}\label{t} Assume \eqref{e:greedy}. Let $\lambda>0$. Then
$$
\mu(\lambda,\nu,d) = 0 \mbox{ if and only if } \lambda \ge \widehat{\lambda}_c(\nu,d).
$$
\end{theorem}

Since the set $\{\lambda>0 :  \lim_{r \rightarrow \infty }P(S(r) \as S(2r)) = 0\}$ is open (see Appendix \ref{s:explicitons}), Theorem \ref{t} is in fact equivalent to the following proposition, that we actually prove in the next sections.
\begin{prop}
\label{prop}
Assume \eqref{e:greedy}. Let $\lambda>0$. Then
$$
\mu(\lambda,\nu,d) > 0 \mbox{ if and only if } \lim_{r \rightarrow \infty }P(S(r) \as S(2r)) = 0.
$$
\end{prop}

Let us define a new threshold by
$$ \lambda_{\mu} = \lambda_{\mu} (\nu, d) = \sup \{ \lambda >0 \,:\, \mu (\lambda, \nu, d) =0  \} \,. $$
Theorem \ref{t} can be reformulated to obtain the following corollary.
\begin{corollary}
\label{cor}
Assume \eqref{e:greedy}. Then
$$\lambda_{\mu} = \widehat \lambda_c.$$
Moreover, 
$$ \mu (\lambda_{\mu} (\nu, d), \nu, d) = 0. $$
\end{corollary}

Theorem \ref{t} is analogous to the result of Kesten \cite{Kesten-saint-flour} (Theorem 6.1) 
in the framework of Bernoulli percolation and first passage percolation on $\Z^d$.
The proof of Kesten can be adapted in our setting in the case of bounded radii.
In the general case, some further arguments are needed.

\bigskip

In \cite{GM-deijfen}, the following result was implicitly proved:
if \eqref{e:greedy} holds, then $\mu(\lambda,\nu,d)$ is positive for small enough $\lambda>0$.
Theorem \ref{t} is therefore a strengthening of this result.

\paragraph{Remark.}
We did not work out the details but we think that the result and the proof can be extended to the following framework.
Replace the random Euclidean balls centered at each point of the Poisson point process by independent copies of a given random subset $S$.
Assume that $S$ is regular enough and assume the existence of a deterministic constant $K$ such that
\begin{itemize}
\item With probability one, there exists $r>0$ such that $B(0,r) \subset S \subset B(0,Kr)$.
\item With probability one, for any $x,y \in \overline{S}$, there exists a polygonal path $\pi$ in $S$, with the possible exception of its first and last point, such that $\ell(\pi) \le K \mbox{diameter}(S)$.
\end{itemize}
The integrability conditions on the law of the radii are then replaced by integrability conditions on the diameter of $S$.

\bigskip

\paragraph{Acknowledgements.} The authors would like to thank an anonymous referee for its many valuable comments that helped to improve and clarify the article in general, and especially to greatly simplify the proof of Lemma \ref{l:geodesic}.

\section{Proof of $\mu>0 \Rightarrow  \lim_{r \rightarrow \infty }P(S(r) \as S(2r)) = 0$} 

Let $\lambda>0$. Let us first prove that
\begin{equation}\label{e:consequence_shape}
\lim_{r \rightarrow \infty}\frac {T(S(r),S(2r))}r = \mu \mbox{ a.s. } 
\end{equation}
Any path from $0$ to $S(2r)$ can be seen as the concatenation of a first path from $0$ to $S(r)$ and a second path from $S(r)$ to $S(2r)$.
Taking infimums, we get
\begin{equation}\label{e:vmc1}
T(0,S(r)) + T(S(r),S(2r)) \le T(0,S(2r)).
\end{equation}
On the other hand, for any $x$ in $S(r)$ we have
\begin{eqnarray*}
T(0,S(2r)) 
 & \le  & T(0,x) + T(x,S(2r)) \\
 & \le & \left(\sup_{x' \in S(r)} T(0,x')\right) + T(x,S(2r)).
\end{eqnarray*}
Taking the infimum in $x$, we now get
\begin{equation}\label{e:vmc2}
T(0,S(2r))  \le \left(\sup_{x' \in S(r)} T(0,x') \right) + T(S(r),S(2r)).
\end{equation}
From \eqref{e:vmc1} and \eqref{e:vmc2} we get
$$
T(0,S(2r)) -\left(\sup_{x' \in S(r)} T(0,x') \right)   \le T(S(r),S(2r)) \le T(0,S(2r)) - T(0,S(r)).
$$
By Theorem \ref{t:forme} we then deduce \eqref{e:consequence_shape}.

\bigskip

Now assume $\mu(\lambda,\nu,d) >0$. Let us prove $\lim_{r \rightarrow \infty }P(S(r) \as S(2r)) = 0$.
For all $r>0$,
$$
P\big(S(r) \as S(2r)\big) \le P\left(\frac{T(S(r),S(2r))}r=0\right).
$$
But this tends to $0$ as $r$ tends to infinity thanks to \eqref{e:consequence_shape} and the assumption $\mu>0$.

\section{Proof of $ \lim_{r \rightarrow \infty }P(S(r) \as S(2r)) = 0 \Rightarrow \mu>0$}

Let $\lambda>0$ such that $ \lim_{r \rightarrow \infty }P(S(r) \as S(2r)) = 0$. In this section, we prove that this implies $\mu(\lambda,\nu,d) > 0$.
Let us give the plan of the proof. First fix a large enough $\rho$. Since $\lim_{r \rightarrow \infty }P(S(r) \as S(2r)) = 0$,
we can suppose that with high probability the time needed to cross an annulus of the form $B(x,2\rho) \setminus B(x,\rho)$
is bigger than some positive constant $\delta$. We consider now a very large $r$, a lot bigger than $\rho$. We consider a nice geodesic $\pi$ from $0$ to the Euclidean sphere $S(r)$. Say, for concreteness, that a random ball of the Boolean model is large if its radius is larger that $10\rho$.

First, imagine that there exists no large random balls.
The idea is then to discretize the geodesic $\pi$ in a convenient way at scale $\rho$:
the traveler, when moving along the geodesic, crosses several annuli of the form $B(x,2\rho) \setminus B(x,\rho)$.
Let us say that such an annulus is good if the time needed to cross the annulus is at least $\delta$.
We have
$$
\tau(\pi) \ge \delta \# \{\mbox{good annuli crossed}\}.
$$
By our choice of $\rho$ and $\delta$, each given annulus is good with high probability.
Moreover, as there are no large balls, what occurs in far enough annuli is independent.
From these observations one can prove
$$
\tau(\pi) \ge \frac{C\delta r}\rho
$$
for some constant $C$ with high probability and the positivity of $\mu$ follows.

The difficulty is to take care of large balls.
Let $\tau_\rho(\cdot)$ denote the time needed to travel along a path when we throw away large balls.
Let us define the notion of good annulus as before using travel times $\tau_\rho$.
The time needed to cross a good annulus is at least $\delta$ unless there is a large ball with touches the annulus.
If the time needed to cross a good annulus is smaller than $\delta$, then some large ball touches the annulus.
In that case, we say that the annulus is disturbed by the large ball.
We get
$$
\tau(\pi) \ge \delta \# \{\mbox{good annuli crossed}\} - \delta \# \{\mbox{disturbed annuli}\}.
$$
We take care of the first term of the right-hand side as before.
We handle the second term - the perturbative term -  by relating it to the greedy paths models.
We show that, for $\rho$ large enough, the perturbative term is smaller than the first term and the positivity of $\mu$ follows.

\subsection{Greedy paths}

Let $\pi=(x_0,\dots,x_k)$ be a path. Recall that a path is a family of distinct points of $\R^d$. Set
$$
r(\pi)=\sum_{1 \le i \le k}  r(x_i)
$$
(recall that $r(x)$ is the radius of the ball centered at $x$, as defined at the end of the paragraph about the Boolean model) and
$$
S = S(\xi) = \sup_{\pi} \frac{r(\pi)}{\ell(\pi)}
$$
where the supremum runs over all paths such that $x_0=0$ and $k \ge 1$.
Note that the intensity measure of the underlying Poisson point process $\xi$ is $\lambda |\cdot| \otimes \nu = |\cdot| \otimes \lambda \nu$.

We can consider the greedy path model when the intensity measure is $|\cdot| \otimes m$ where $m$ is a given finite measure on $(0,+\infty)$.
In that case, we write $r_m(\pi)$ and $S_m$.

\begin{theorem}[\cite{GM-deijfen}] \label{t:greedy} Let $m$ be a finite measure on $(0,+\infty)$. There exists a constant $C=C(d)$ such that
$$
E\left(S_m\right) \le C \int_{(0,\infty)} m\big((r,+\infty)\big)^{1/d} dr.
$$
\end{theorem}
This is a consequence of (11) in \cite{GM-deijfen} and Lemma 2.1 in the same article.
Note that the results requires the assumption $d \ge 2$.
The result is the analogue in the continuous setting of a result by Martin \cite{Martin-greedy} in the discrete setting.

\subsection{Geodesics}

\begin{lemma} \label{l:geodesic} For any $r>0$ there exists a path $\pi=(x_0,\dots,x_k)$ such that
 \begin{enumerate}
  \item[(i).] $x_0=0$, $\|x_k\|=r$ and, for all $i$, $\|x_i\| \le r$. 
  \item[(ii).] $\tau(\pi)=T(r)$.
  \item[(iii).] For any ball $B=B(c,r(c))$ of the Boolean model, $\pi^{-1}(B)$ (here we see $\pi$ as a curve parametrized by arc-length) is contained in an interval of length at most $3r(c)$.
 \end{enumerate}
\end{lemma}

\paragraph{Proof.} Let us call geodesic a path $\pi$ satisfying (i) and (ii). 
We work on the full probability event "the number of random balls which touch $\overline{B(0,r)}$ is finite". 
Let $\Sigma(r)$ denote the union of all random balls which touch $\overline{B(0,r)}$.
Let $V$ be the set of connected components of $\Sigma(r)$.
Define $V'$ as the union of $V$ and $\{\{0\},S(r)\}$.
The set $V'$ is finite.
We consider the complete graph whose vertices set is $V'$.
The length of an edge between any $C,C'$ in $V'$ is defined as the Euclidean distance between $C$ and $C'$.
We consider the natural associated geodesic distance $d$ on the graph.
Let us check
\begin{equation}\label{e:geodesique-discretisation}
T(r) = d(\{0\},S(r))
\end{equation}
and the existence of a geodesic.

\begin{figure}[h!]
\centering
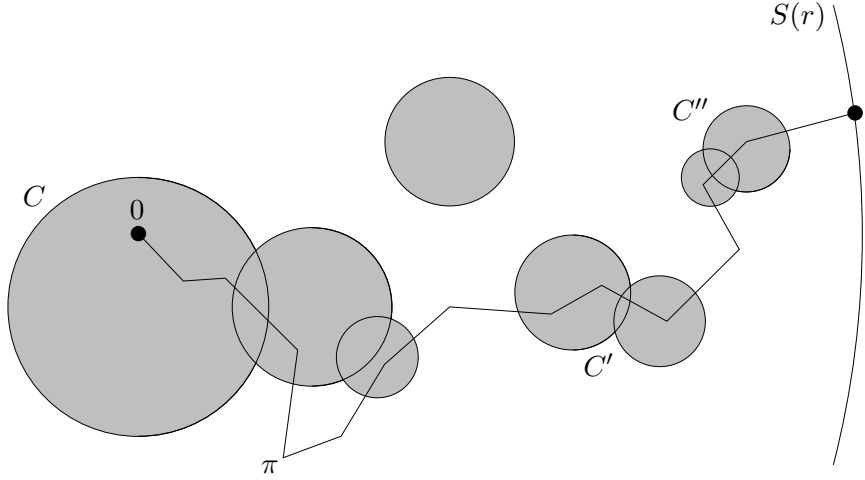
\caption{A path $\pi$ from $0$ to $S(r)$, and the corresponding sequence $(C(0), C(1), C(2), C(3), C(4), C(5)) = (\{0\}, C, C, C', C'', S(r))$.}
\label{fig1}
\end{figure}

We first prove the inequality $T(r) \ge d(\{0\},S(r))$.
Let $\pi$ be a path from $0$ to $S(r)$ whose image is contained in $\overline{B(0,r)}$.
We see $\pi$ as a curve $[0,\ell(\pi)] \to \R^d$ parametrized by arc-length.
Let $(C(1),\dots,C(n-1))$ be the finite sequence 
\footnote{The definition makes sense because of the following facts.
\begin{itemize}
 \item The set $V$ is finite.
 \item The sets $\pi^{-1}(C), C \in V$ are disjoint.
 \item For each $C \in V$, the set $\pi^{-1}(C)$ is the union of finite number of intervals.
\end{itemize}}
of elements of $V$ successively visited by $\pi$. We refer to Figure \ref{fig1}.

Set $C(0)=\{0\}$ and $C(n)=S(r)$.
Thus, $(C(0),\dots,C(n))$ is a sequence of elements of $V'$.
For any $i \in \{0,\dots,n-1\}$, some part (possibly empty) of $\pi$ goes from $C(i)$ to $C(i+1)$ without touching $\Sigma$.
Here again we refer to Figure \ref{fig1}.
The travel time of this part of $\pi$ is at least $d(C(i),C(i+1))$ and
$$
\tau(\pi) \ge \sum_{i=0}^{n-1} d(C(i),C(i+1) \ge d(\{0\},S(r)).
$$
Therefore
$$
T(r) \ge d(\{0\},S(r)).
$$

\begin{figure}[h!]
\centering
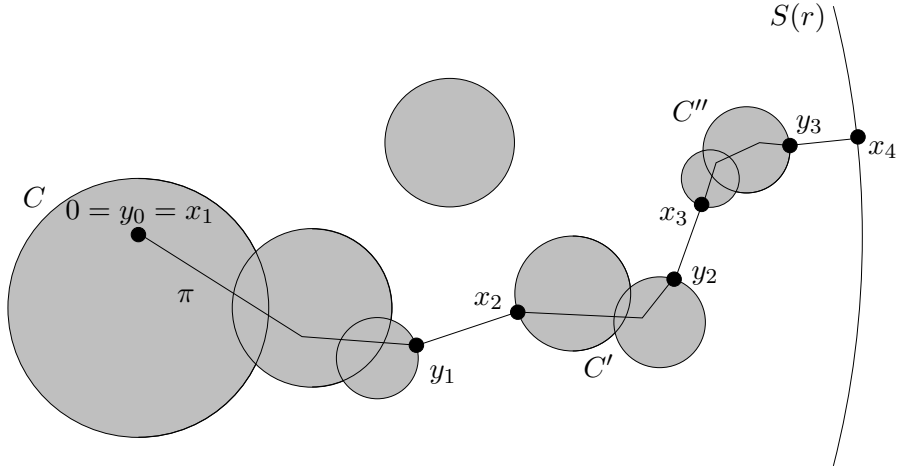
\caption{The sequence $(C(0), C(1), C(2), C(3), C(4)) = (\{0\}, C, C', C'', S(r))$ of vertices of $V'$, the corresponding points $(y_i, x_{i+1})_{0\leq i \leq n-1}$ and a geodesic $\pi$.}
\label{fig2}
\end{figure}

We now prove the reverse inequality by constructing a geodesic for $T(r)$. We can construct such a geodesic as follows.
Let $(C_0=\{0\},C(1),\dots,C(n)=S(r))$ be a sequence of distinct vertices of $V'$ such that
$$
d(\{0\},S(r)) = \sum_{i=0}^{n-1} d(C(i),C(i+1))
$$
(see Figure \ref{fig2}).
For any $i \in \{0,\dots,n-1\}$, let $y_i \in \overline{C(i)}$ and $x_{i+1} \in \overline{C(i+1)}$ be such that $d(C(i),C(i+1))=\|x_{i+1}-y_i\|$ (see Figure \ref{fig2}).
Then, concatenate the following polygonal paths : 
the segment from $y_0=0$ to $x_1$ ($0=x_1$ can occur), 
a path from $x_1$ to $y_1$ in $C_1$ (with the exception of the initial and final points), 
the segment from $y_1$ to $x_2$ (here again $y_1=x_2$ can occur),
and so on until the point $x_n \in S(r)$.
This is a path  $\pi$ from $0$ to $S(r)$ such that 
$$
\tau(\pi) \leq \sum_{i=0}^{n-1} \|x_{i+1}-y_i\| =  \sum_{i=0}^{n-1} d(C(i),C(i+1)) = d(\{0\},S(r)).
$$
Therefore  \eqref{e:geodesique-discretisation} holds and $\pi$ (or $\pi$ stopped at its first encounter with $S(r)$) is a geodesic.

\medskip

Set
$$
\ell = \inf\{ \ell(\pi), \pi \mbox{ is a geodesic}\}
$$
and let $\pi$ be a geodesic such that
$$
\ell(\pi) \le \ell + r_0
$$
where $r_0$ is the minimum of the radii of all balls of the Boolean model which touch $\overline{B(0,r)}$.
Let $B=B(c,r(c))$ be a ball of the Boolean model.
Assume that $\pi^{-1}(B)$ is not empty, otherwise there is nothing to prove.
Then $B$ touches $\overline{B(0,r)}$ and therefore $r(c) \ge r_0$.
Let $s=\inf \{\pi^{-1}(B)\}$ and $t=\sup \{\pi^{-1}(B) \}$.
As $\pi(s)$ and $\pi(t)$ belong to $\overline{B}$, the travel time of the segment from $\pi(s)$ to $\pi(t)$ equals $0$.
Modify $\pi$ by replacing the part of $\pi$ between $s$ and $t$ by the segment from $\pi(s)$ to $\pi(t)$.
Let $\widetilde\pi$ denote the new path.
We also see it as parametrized by arc-length.
By the previous remark and as $\pi$ is a geodesic, $\widetilde\pi$ is also a geodesic and therefore
$$
\ell \le \ell(\widetilde{\pi}).
$$
But,
$$
\ell(\widetilde{\pi}) \le \ell(\pi)-(t-s)+\|\pi(t)-\pi(s)\| \le \ell+r_0-(t-s) + 2r(c) \le \ell + 3r(c) - (t-s). 
$$
The two previous inequalities yield
$$
t-s \le 3r(c).
$$
The lemma follows. \finpreuve

\paragraph{Remarks.} 
\begin{itemize}
 \item One can prove the existence of a geodesic $\pi$ such that $\ell(\pi)$ is minimal among all geodesics. Such a geodesic $\pi$ satisfies
 \begin{itemize}
 \item[(iii').] For any ball $B=B(c,r(c))$ of the Boolean model, $\pi^{-1}(B)$ is contained in an interval of length at most $2r(c)$.
\end{itemize}
 \item There exists geodesics which (with the exception of the first and last point) goes from one center of a random ball to one other. But in some cases, no such geodesics satisfies (iii'). 
 \item We decided to state and prove Lemma \ref{l:geodesic} this way for two reasons.
First, it is sufficient for our purpose and it is easy to prove.
Second, the proof can easily be adapted to a more general case where, for example, random balls are replaced by random shapes $S$ fulfilling some regularity properties together with the following property:
for any $x,y \in \overline{S}$, there exists a polygonal path $\pi$ from $x$ to $y$ inside $S$ (with the possible exception of the first and last point) 
such that $\ell(\pi)$ is at most $K$ times the diameter of $S$ for a uniform constant $K$.
\end{itemize}

\subsection{Constants}
\label{sec:constant}

The aim of this section is to define the following constants : $\eta, K, \kappa, \epsilon, \rho$ and $\delta$. The reader can skip this section at first reading.

\medskip

Fix $\eta=\eta(d)>0$ such that
\begin{equation}\label{e:definition:eta}
\R^d \subset \bigcup_{i \in \Z^d} B(\eta i,1).
\end{equation}
Let $K=K(d)$ be the cardinality of $B(0,12 \eta^{-1}) \cap \Z^d$.
Fix $\kappa=\kappa(d)>0$ large enough and $\epsilon=\epsilon(d)>0$ small enough such that
\begin{equation}\label{e:definition:kappaeps}
K\exp(-\kappa/2)+K\exp(\kappa/2)\epsilon \le 1/2.
\end{equation}
Fix $\rho=\rho(\lambda,\nu,d)>0$ large enough such that
$$
P(S(\rho) \as S(2\rho)) \le \epsilon/2
$$
and
\begin{equation}\label{e:A-greedy}
C \int_{(0,+\infty)} \lambda^{1/d}\nu^\rho \big((r,+\infty)\big)^{1/d} dr \le 4000^{-1}
\end{equation}
where $C$ is the constant which appears in Theorem \ref{t:greedy} and where the measure $\nu^\rho$ is defined by $\nu^\rho(\cdot)=\nu(\cdot \cap [\rho,+\infty)])$. 
Since the following inclusion holds almost surely\footnote{Here is one way to prove the inclusion. There exists a geodesic $\pi$ such that $T(S(\rho), S(2\rho)) = \tau (\pi)$. 
This can be shown for example in the same way as Lemma \ref{l:geodesic}. The almost sure inclusion  $\{S(\rho) \as S(2\rho)\}^c \subset \{T(S(\rho),S(2\rho)) > 0\}$ follows from the existence of the geodesic and the fact that almost surely there does not exist two random balls that intersect $B(2\rho)$ and are tangent.}
$$\{S(\rho) \as S(2\rho)\}^c \subset \{T(S(\rho),S(2\rho)) > 0\}.$$
we have
$$
P(T(S(\rho),S(2\rho)) > 0) \ge 1-\epsilon/2.
$$
Therefore we can fix $\delta=\delta(\lambda,\nu,d)>0$ such that
\begin{equation}\label{e:delta}
P( T(S(\rho), S(2\rho)) \ge \delta) \ge 1-\epsilon.
\end{equation}

\subsection{Skeletons of the geodesic.}

Let $r \ge 20 \rho$.
Let $\pi$ be the geodesic given by Lemma \ref{l:geodesic}.
We see $\pi$ as a curve $[0,\ell(\pi)] \to \R^d$ parametrized by arc length.

\begin{lemma} \label{l:skeleton} Set 
\begin{equation}\label{e:k}
k=\lceil r/(20\rho) \rceil.
\end{equation}
There exists a sequence $0=t(0) \le t(1) \le \cdots \le t(k) \le \ell(\pi)$ such that
\begin{enumerate}
 \item For any distinct $j,j' \in \{0,\dots,k\}$,  $\|\pi(t(j))-\pi(t(j'))\| \ge 10\rho$.
 \item For any  $j \in \{1,\dots,k\}$, $\|\pi(t(j))-\pi(t(j-1))\| \le 10\rho$.
\end{enumerate}
\end{lemma}

\paragraph{Proof. }
We build a sequence $0=t(0) \le t(1) \cdots \le t(k) \le \ell(\pi)$ as follows.
Let $t(0)=0$ and $k=0$.
We proceed by induction.
At each step we consider the set
$$
\left\{t \in [0,\ell(\pi)] : \pi(t) \in \bigcup_{j \le k} \overline{B}\big(\pi(t(j)),10\rho\big) \right\}.
$$
\begin{itemize}
\item If $\ell(\pi)$ belongs to the set, then the construction is over.
\item Otherwise, we define $t(k+1)$ as the maximum of this non empty compact set, we increase by one $k$ and the construction goes on.
\end{itemize}
By throwing away, if needed, the last elements of our sequence, we get a new sequence which fulfills the properties stated in the Lemma.
Let us give some details.
\begin{itemize}
\item The sequence is clearly non-decreasing. It is  actually increasing. Indeed, let $j' \in \{1,\dots,k\}$. By construction, the interval  $]t(j'),\ell(\pi)]$ is non-empty.
But if $t(j')=t(j'-1)$, then there would exists $t$ in the previous interval such that $\pi(t) \in \overline{B}(\pi(t(j'-1)),10\rho)$. This would contradict the definition of $t(j')$.
\item Property 1 of the lemma holds. Indeed, let $j,j' \in \{0,\dots,k\}$ be such that $j < j'$. By construction, the interval  $]t(j'),\ell(\pi)]$ is non-empty and for any $t$ in this interval, $\pi(t) \not\in\overline{B}(\pi(t(j)),10\rho)$.
The result follows by continuity of $\pi$.
\item Property 2 of the lemma holds. Indeed, let $j' \in \{1,\dots,k\}$. By construction, there exists $j \in \{0,\dots,j'-1\}$ such that $\pi(t(j')) \in \overline{B}(\pi(t(j)),10\rho)$.
If $j \le j'-2$,  by construction of $t(j'-1)$, one would have $t(j'-1) \ge t(j')$ which is not true. Therefore $j=j'-1$.
\end{itemize}
Thanks to the stopping criterion, $\pi(\ell(\pi)) \in \overline{B}(\pi(t(j)),10\rho)$ for some $j \le k$.
As $\|\pi(\ell(\pi))\|=r$, Property 2 of the lemma yields
$$
(k+1)10\rho \ge (j+1)10\rho \ge r.
$$
As $r \ge 20 \rho$ we get $k \ge 1$ and then $20k \rho \ge r$.
Thus, by throwing away, if needed, the last elements of our sequence, we get a new sequence which fulfills the properties stated in the Lemma.\finpreuve

\subsection{About disturbant balls}

For all $j \in \{0,\dots,k\}$, we set
$$
\Pi(j)=\pi\big([0,\ell(\pi)]\big) \cap B(t(j),3\rho)
$$
where $\pi$ is still the geodesic given by Lemma \ref{l:geodesic} seen as a curve parametrized by arc length.
Let us say that a random ball $B(c,r(c))$ of the Boolean model disturbs $\Pi(j)$ 
if $\Pi(j) \cap B(c,r(c)) \neq \emptyset$ and if $r(c) \ge \rho$.
If no random ball disturbs $\Pi(j)$, then the travel time of the geodesic inside $B(t(j),3\rho)$ does not depend on large balls.
Set 
$$
D(j) = \{\mbox{centers of random balls which disturb }\Pi(j)\}
$$
and
$$
D = \bigcup_{j \in \{0,\dots, k\}} D(j).
$$
We will use the following upper bound on the number of $\Pi(j)$ which are disturbed by a given large ball of the Boolean model.

\begin{lemma} \label{l:disturb} Let $B(c,r(c))$  be a ball of the Boolean model such that $r(c) \ge \rho$. 
Then
$$
\card(\{j \in \{0,\dots,k\} : \Pi(j) \cap B(c,r(c)) \neq \emptyset \}) \le \frac{2r(c)}\rho.
$$
\end{lemma}

\paragraph{Proof.} The idea is that the geodesic restricted to $B(c,r(c))$ is contained in a path of length at most $3r(c)$ 
and that points in different $\Pi(j)$ are at least at distance $4\rho$ from each other. Thanks to these observations we will prove that
$$
\card(\{j \in \{0,\dots,k\} : \Pi(j) \cap B(c,r(c)) \neq \emptyset \}) \le \frac{3r(c)}{4\rho}+1
$$
and this inequality enables us to conclude.

Let us give a detailed proof.
For all $j \in \{0,\dots,k\}$, we set
$$
U(j) = \pi^{-1}\big(B(t(j),3\rho)\big).
$$
Note that $\Pi(j)=\pi(U(j))$.
Set 
$$
J = \{j \in \{0,\dots,k\} : \pi(U(j)) \cap B(c,r(c)) \neq \emptyset\}.
$$
We aim at proving $\card(J) \le 2r(c)/\rho$.

For all $j \in J$, fix $t'(j) \in U(j)$ such that $\pi(t'(j)) \in B(c,r(c))$.
As $t'(j) \in U(j)$, we have $\pi(t'(j)) \in B(t(j), 3\rho)$.
For any distinct $j,j' \in J$, by Lemma \ref{l:geodesic}, we thus have 
$$
|\pi(t'(j))-\pi(t'(j'))| \ge |\pi(t(j))-\pi(t(j'))|-6\rho \ge 10\rho - 6\rho=4\rho.
$$
As $\pi$ is parametrized by arc length, we then get 
\begin{equation}\label{e:quivachercherFlorence}
|t'(j)-t'(j')| \ge 4\rho.
\end{equation}
By Lemma \ref{l:geodesic}, $\pi^{-1}(B(c,r(c)))$ is contained in an interval of length at most $3r(c)$.
Ordering the $t'(j), j \in J$, using the fact that each such $t'(j)$ belongs to $\pi^{-1}(B(c,r(c)))$ and using \eqref{e:quivachercherFlorence} we thus get
$$
3 r(c) \ge \max_j t'(j) - \min_j t'(j) \ge 4\rho(\card(J)-1).
$$
Therefore
$$
\card(J) \le \frac{3r(c)}{4\rho}+1 \le \frac{2r(c)}\rho
$$
as $r(c) \ge \rho$. \finpreuve

\subsection{A second path}

We define a new path $\widetilde\pi$. The definition of $\widetilde{\pi}$ is in a sense artificial, since it is built to enable the use of the results on greedy paths.
\begin{itemize}
 \item It starts from $\pi(t(0)) = 0$ and visits each point of $D(0)$ (see the subsection about disturbing balls), if any.
 \item Then it goes to $\pi(t(1))$ and visits each point of $D(1)$ it has not visited yet, if any.
 \item Then it goes to $\pi(t(2))$ and visits each point of $D(2)$ it has not visited yet, if any.
 \item $\dots$
 \item Then it goes to $\pi(t(k))$ and visits each point of $D(k)$ it has not visited yet, if any.
\end{itemize}
In particular, $\widetilde\pi$ visits all points $\pi(t(j))$, $0 \le j \le k$. By Lemma \ref{l:skeleton} we get
$$
\ell(\widetilde{\pi}) \ge 10\rho k.
$$
As $k \ge r/(20\rho)$, we get
\begin{equation}\label{e:lpitildelower}
\ell(\widetilde{\pi}) \ge r/2.
\end{equation}
If a random ball $B(c,r(c))$ disturbs a set $\Pi(j)$, then
\begin{equation}\label{e:coucouMarie}
\|c-\pi(t(j))\| \le 3\rho+r(c) \le 4r(c).
\end{equation}
We can easily give an upper bound on the length of $\widetilde{\pi}$ by considering the longer sequence in which in the definition 
"visits each point of $D(j)$ it has not visited yet" is replaced by
"goes back and forth between $\pi(t(j))$ and points of $D(j)$ it has not visited yet".
Using \eqref{e:coucouMarie} and Lemma \ref{l:skeleton}, we get
\begin{equation}\label{e:lpitildeupper}
\ell(\widetilde{\pi}) \le 10\rho k+2\sum_{c \in D} 4r(c) = 10\rho k+8\sum_{c \in D} r(c). 
\end{equation}

\subsection{Lower bound on $T(r)$ (1/2): large balls influence can be controlled by greedy paths}

Set
$$
\Sigma_- = \bigcup_{(c,r) \in \xi : r < \rho} B(c,r).
$$
In other words, we throw away all balls of radius larger than or equal to $\rho$.
Let $T_-$ be defined from $\Sigma_-$ in the same way as $T$ is defined from $\Sigma$.
In particular $T_- \ge T$. We list a few facts.
\begin{itemize}
 \item For any $j \in \{0,\dots,k-1\}$, the path $\pi$ contains a path from $\pi(t(j))$ to $S(\pi(t(j)), 3\rho)$.
This is a consequence of the first part of Lemma \ref{l:skeleton} which yields, for any such $j$, $\pi(t(j+1)) \not\in B(\pi(t(j)), 3\rho)$.
 \item The balls $B(\pi(t(j)), 3\rho)$, $j \in \{0,\dots,k-1\}$ are disjoint. This is again a consequence of the first part of Lemma \ref{l:skeleton}.
 Therefore the travel length of $\pi$ is at least the sum of the $T(\pi(t(j)), S(\pi(t(j)),3\rho)), j \in \{0,\dots,k-1\}$.
 \item For any $j \in \{0,\dots,k-1\}$, if $D(j)$ is empty, then the travel time of $\pi$ inside $B( \pi(t(j)),3\rho)$ does not change if we throw away random balls of radii at least $\rho$.
\end{itemize}
As a consequence,
$$
T(r) \ge \sum_{j=0}^{k-1} T_-(\pi(t(j)),S(\pi(t(j)),3\rho)) 1_{D(j)=\emptyset}.
$$
Therefore
\begin{eqnarray*}
T(r)
 & \ge & \delta \sum_{j=0}^{k-1} 1_{T_-(\pi(t(j)),S(\pi(t(j)),3\rho))\ge \delta} 1_{D(j)=\emptyset} \\
 & \ge & \delta \sum_{j=0}^{k-1} 1_{T_-(\pi(t(j)),S(\pi(t(j)),3\rho))\ge \delta} - \delta \sum_{j=0}^{k-1} 1_{D(j) \neq \emptyset}.
\end{eqnarray*}
  By Lemma \ref{l:disturb}, a given random ball $B(c,r(c))$ such that $r(c) \ge \rho$ disturbs at most $2r(c)/\rho$ sets $\Pi(j)$.
Therefore
$$
T(r) \ge  \delta \sum_{j=0}^{k-1} 1_{T_-(\pi(t(j)),S(\pi(t(j)),3\rho)) \ge \delta} - \delta \sum_{c \in D} \frac{2r(c)}\rho
$$
where we recall that $D$ is the union of $D(j), j \in \{0,\dots, k\}$.
Recall that any point of $D$ is an element of the path $\widetilde{\pi}$ seen as a sequence.
Note also that any point of $D$ is the center of a ball of radius at least $\rho$.
Therefore, using the notations of the greedy paths model,
\begin{equation}\label{e:liengreedy}
\sum_{c \in D} r(c) \le r_{\lambda \nu^\rho}(\widetilde\pi)
\end{equation}
where
$$
\nu^\rho = \nu(\cdot \cap [\rho,+\infty))
$$
and where we use the natural coupling between the Poisson point process which defines the Boolean model 
and the Poisson point process which defines the above greedy paths model:
we just throw away all point $(c,r)$ such that $r<\rho$.
Therefore 
$$
T(r) \ge  \delta \sum_{j=0}^{k-1} 1_{T_-(\pi(t(j)),S(\pi(t(j)),3\rho))\ge \delta} - \frac{2\delta}\rho r_{\lambda \nu^\rho}(\widetilde{\pi}).
$$
At this point we would like to divide $T(r)$ by $r$, the first sum on the right-hand side by $k$ and the second sum by $\ell(\widetilde{\pi})$. 
By definition of $k$ (see \eqref{e:k}), $r$ and $k$ are closely related. 
The links between $\ell(\widetilde{\pi})$ and $r$ (or $k$) is less clear.
Therefore we use the following trick.
$$
T(r) \ge  \delta \sum_{j=0}^{k-1} 1_{T_-(\pi(t(j)),S(\pi(t(j)),3\rho)) \ge \delta}+\frac{16\delta}\rho r_{\lambda \nu^\rho}(\widetilde{\pi})  
- \frac{18\delta}\rho r_{\lambda \nu^\rho}(\widetilde{\pi}).
$$
By \eqref{e:lpitildelower} we get
$$
\frac{T(r)}r \ge \frac{T(r)}{2\ell(\widetilde{\pi})}
$$
and then
\begin{eqnarray*}
\frac{T(r)}r  
& \ge  & \delta \left( \frac{\sum_{j=0}^{k-1} 1_{T_-(\pi(t(j)),S(\pi(t(j)),3\rho)) \ge \delta} 
+ \frac{16}\rho r_{\lambda \nu^\rho}(\widetilde{\pi})}{2\ell(\widetilde{\pi})} 
- \frac{\frac{18}\rho r_{\lambda \nu^\rho}(\widetilde{\pi})}{2\ell(\widetilde{\pi})} \right)\\
& = & \frac{\delta}{\rho} 
\left(
\frac{\rho \sum_{j=0}^{k-1} 1_{T_-(\pi(t(j)),S(\pi(t(j)),3\rho)) \ge \delta} + 16 r_{\lambda \nu^\rho}(\widetilde{\pi})}{2\ell(\widetilde{\pi})} 
-  \frac{9r_{\lambda \nu^\rho}(\widetilde{\pi})}{\ell(\widetilde{\pi})} 
\right).
\end{eqnarray*}
Using \eqref{e:lpitildeupper} and \eqref{e:liengreedy} we then get
$$
\frac{T(r)}r \ge \frac{\delta}{ \rho} 
\left(
\frac{\rho \sum_{j=0}^{k-1} 1_{T_-(\pi(t(j)),S(\pi(t(j)),3\rho)) \ge \delta} 
+ 16 r_{\lambda \nu^\rho}(\widetilde{\pi})}{20\rho k+16r_{\lambda \nu^\rho}(\widetilde{\pi})}
- \frac{ 9 r_{\lambda \nu^\rho}(\widetilde{\pi})}{\ell(\widetilde{\pi})}
\right).
$$
Using
$$
\rho \sum_{j=0}^{k-1} 1_{T_-(\pi(t(j)),S(\pi(t(j)),3\rho)) \ge \delta} \le 20\rho k
$$
we obtain
\begin{eqnarray*}
\frac{T(r)}r 
 & \ge & 
\frac{\delta}{ \rho} 
\left(
\frac{\rho \sum_{j=0}^{k-1}1_{T_-(\pi(t(j)),S(\pi(t(j)),3\rho)) \ge \delta} 
}{20\rho k}
- \frac{ 9 r_{\lambda \nu^\rho}(\widetilde{\pi})}{ \ell(\widetilde{\pi})}
\right) \\
 & \ge & \frac{\delta}{ \rho} 
\left(
\frac{ \sum_{j=0}^{k-1} 1_{T_-(\pi(t(j)),S(\pi(t(j)),3\rho)) \ge \delta} 
}{20k}
- \frac{10  r_{\lambda \nu^\rho}(\widetilde{\pi})}{\ell(\widetilde{\pi})}
\right).
\end{eqnarray*}
Using notations for greedy paths, we thus get
\begin{equation}\label{e:onreecrittout}
\frac{T(r)}r 
\ge 
\frac{\delta}{ \rho} 
\left(
 \frac{ \sum_{j=0}^{k-1} 1_{T_-(\pi(t(j)),S(\pi(t(j)),3\rho)) \ge \delta}  }{20k}
 - 10S_{\lambda \nu^\rho}
\right).
\end{equation}

\subsection{Lower bound on $T(r)$ (2/2): control of small balls influence and conclusion}

Let ${\mathcal S}$ be the set of sequences $(s(0),\dots,s(k-1))$ of distinct elements of $\eta \rho \Z^d$ such that
\begin{itemize}
 \item $s(0)=0$.
 \item For all distinct $j,j' \in \{0,\dots,k-1\}$,  $\|s(j)-s(j')\| \ge 8 \rho$.
 \item For all $j \in \{1,\dots,k-1\}$, $\|s(j)-s(j-1)\| \le 12 \rho$.
\end{itemize}
We now define a new sequence $s(0),\dots,s(k-1)$ which is a discretization of the sequence $\pi (t(0)),\dots,\pi(t(k-1))$.
We set $s(0)=0$.
For all $j \in \{1,\dots,k-1\}$, we chose  $s(j) \in \eta \rho \Z^d$ such that $\pi(t(j))$ belongs to $B(s(j),\rho)$.
This is possible thanks to \eqref{e:definition:eta}.

\begin{lemma}\label{l:skeleton2} The sequence $s(0), \dots, s(k-1)$ fulfills the following properties.
\begin{enumerate}
 \item For all $j \in \{0,\dots,k-1\}$, $\pi(t(j))$ belongs to $B(s(j) ,\rho)$.
 \item The sequence $s$ belongs to ${\cal S}$.
\end{enumerate}
\end{lemma}
\paragraph{Proof.} The first property holds by definition for any index $j \ge 1$. It also holds for $j=0$ as $\pi(t(0))=\pi(0)=0$.
Let $j,j'$ be two distinct elements of $\{0,\dots,k-1\}$.
By Lemma \ref{l:skeleton}, we have 
$$
\|\pi(t(j))-\pi(t(j'))\| \ge 10\rho.
$$
By definition of $s(j)$ and $s(j')$, we get $\|s(j) - s(j') \| \ge 8 \rho$.
The remaining property is proven in the same way. \finpreuve

\bigskip

For any $i \in \eta \rho \Z^d$, we say that Site $i$ is good if $T_-(S(i,\rho), S(i,2\rho)) \ge \delta$.
For any $j \in \{0,\dots,k-1\}$, we have $\pi(t(j)) \in B(s(j),\rho)$ and $B(s(j), 2\rho) \subset B(\pi(t(j)), 3 \rho)$ thus $T_-(\pi(t(j)),S(\pi(t(j)),3\rho)) \geq T_-(S(s(j),\rho), S(s(j),2\rho)) $. This implies that if Site $s(j)$ is good then $T_-(\pi(t(j)),S(\pi(t(j)),3\rho)) \ge \delta$.
Therefore \eqref{e:onreecrittout} yields
$$
\frac{T(r)}r 
\ge 
\frac{\delta}{ \rho} 
\left(
 \frac{ \sum_{j=0}^{k-1} 1_{s(j) \mbox{ is good}}}{20k}
 - 10S_{\lambda \nu^\rho}
\right).
$$
and then
\begin{equation}\label{e:lucas}
\frac{T(r)}r 
\ge 
\frac{\delta}{ \rho} 
\left(
 \frac{ \inf_{s' \in {\cal S}} \sum_{j=0}^{k-1} 1_{s'(j) \mbox{ is good}}}{20k}
 - 10S_{\lambda \nu^\rho}
\right).
\end{equation}
By stationarity, by the inequality $T_- \ge T$ and by \eqref{e:delta}, a given site is good with probability at least $1-\epsilon$.
Moreover, as we only consider random balls with radius smaller that $\rho$, the state of Site $i$ only depends on balls whose center belongs to $B(i, 3\rho)$.
Therefore, the state of Sites $i$ and $j$ are independent as soon as $\|i-j\| \ge 6$.
The following lemma is an easy consequence of this observation.

\begin{lemma} The following inequality holds.
\begin{equation}\label{e:I}
P\left(\inf_{s' \in {\mathcal S}} \frac{1}{k}\sum_{j=0}^{k-1} 1_{s'(j) \mbox{ is a good site }} \ge \frac 1 2\right) \ge \frac12.
\end{equation}
\end{lemma}

\paragraph{Proof.} 
Let $A$ denote the complement of the event involved in \eqref{e:I}.
On $A$, there exists at least one path $s' \in {\mathcal S}$ such that
\begin{equation}\label{e:champignon}
\frac{1}{k}\sum_{j=0}^{k-1} 1_{s'(j) \mbox{ is a good site }} \le \frac 1 2.
\end{equation}
For any sequence $s' \in {\mathcal S}$, the events $\{s'(j) \mbox{ is a good site}\}, 0 \le j \le k-1$ are independent and each of them occurs with probability at least $1-\epsilon$. 
This is a consequence of the definition of ${\mathcal S}$ and of the remarks above the lemma.
The probability of the event \eqref{e:champignon} is therefore bounded from above by
$$
P\left(\frac{1}k\sum_{j=0}^{k-1} Z_j \le \frac 1 2\right)
$$
where the $Z_j$ are i.i.d.r.v. with Bernoulli distribution of parameter $1-\epsilon$.
Therefore
\begin{eqnarray*}
 P(A) 
  & \le & \card({\mathcal S}) P\left(\sum_{j=0}^{k-1} Z_j \le k/2\right) \\
  & \le & K^{k-1} P\left(\exp \big(-\kappa \sum_{j=0}^{k-1} Z_j \big)\geq \exp (-k/2)\right)\\
  & \le & K^{k-1}\Big(\exp(\kappa/2)\big((1-\epsilon)\exp(-\kappa) + \epsilon\big)\Big)^k \\
  & \le & \big(K\exp(-\kappa/2)+K\exp(\kappa/2)\epsilon\big)^k \\
  & \le & \frac 1 {2^k}
\end{eqnarray*}
by definition of $K$ and ${\cal S}$, by the choices of $\kappa$ and $\epsilon$ we made (see Section \ref{sec:constant}) and by \eqref{e:definition:kappaeps}.
As $k \ge 1$, the result follows. \finpreuve

\bigskip

By Theorem \ref{t:greedy} and by \eqref{e:A-greedy} we get
$$
P\left(S_{\lambda \nu^\rho} \ge 1/1000\right) \le 1000 E\left(S_{\lambda\nu^\rho}\right) 
\le 1000 C \int_{(0,\infty)} \Big(\lambda \nu^\rho((r,+\infty))\Big)^{1/d} dr \le \frac 1 4.
$$
Therefore, by \eqref{e:I},
$$
P\left(S_{\lambda \nu^\rho} < 1/1000 \mbox{ and } \frac{\inf_{s' \in {\mathcal S}} \sum_{j=0}^{k-1} 1_{s'(j) \mbox{ is a good site}} }{k}  \ge \frac 12 \right) \ge \frac 14
$$
and by \eqref{e:lucas} we obtain
$$
P\left(\frac{T(r)}r \ge \frac{\delta}{80\rho}\right) \ge \frac 14.
$$
As $T(r)/r$ converges in probability to $\mu$, we get $\mu>0$.


\appendix
\appendixpage

\section{Openness of $\{\lambda >0 : \lim_{r \rightarrow \infty}P(S(r) \as S(2r)) = 0\}$ and positivity of $\widehat\lambda_c$}

\label{s:explicitons}

The aim of this section is to provide a proof of the following result. Recall that we assume \eqref{e:momentd}.

\begin{theorem} \label{t:explicitons} The set
 $$
 \{\lambda>0 : \lim_{r\rightarrow \infty} P(S(r) \as S(2r)) = 0\}
 $$
 is open and non-empty. In particular, $\widehat\lambda_c$ is positive.
\end{theorem}

The positivity of $\widehat\lambda_c$ is implicit in \cite{G-perco-boolean-model}.
The openness is a simple consequence of intermediate results in \cite{G-perco-boolean-model}.
Both results are also consequences of Theorems 2.7 and 2.8 in \cite{G-perco-generale} which deal with a more general framework.

We choose to give a proof using intermediate results in \cite{G-perco-boolean-model}.
There is essentially no novelty in this section.

\bigskip

Let us recall some notation from \cite{G-perco-boolean-model}.
Let $\alpha>0$.  
\begin{itemize}
 \item $\Sigma(B(0,\alpha))$ is the union of random balls of the Boolean model with centers in $B(0,\alpha)$.
 \item $G(0,\alpha)$ is the event "there exists a path from $S(\alpha)$ to $S(8\alpha)$ in $\Sigma(B(0,10\alpha))$".
 \item $H(\alpha)$ is the event "there exists a random ball of the Boolean model with which touches $B(0,9\alpha)$ and whose center is outside $B(0,10\alpha)$".
 \item $\Pi(\alpha)=P(G(0,\alpha))$.
\end{itemize}

\bigskip

The article \cite{G-perco-boolean-model} focus on the property $\lim_{\alpha \rightarrow \infty} \Pi(\alpha) = 0$ while in this article we focus on the property 
$\lim_{\alpha \rightarrow \infty}P(S(\alpha) \as S(2\alpha)) = 0$.
This is only a matter of taste, as shown by the first part of the following proposition.
Set 
\begin{equation} \label{e:application-epsilon}
 \epsilon(\alpha)=\int_{[\alpha,+\infty)} r^d \nu(dr).
\end{equation}
Note that $\lim_{\alpha \rightarrow \infty}\epsilon(\alpha) =0$.

\begin{prop} \label{l:prop} There exists a constant $K=K(d)$ such that, for any $\alpha > 0$,
\begin{eqnarray} 
\Pi(\alpha) & \le & P(S(\alpha) \as S(2\alpha)) \le K\Pi(\alpha/10) + \lambda K \epsilon(\alpha/10), \label{e:olivier} \\
\Pi(10\alpha) & \le & K\Pi(\alpha)^2+\lambda K \epsilon(\alpha), \label{e:florence} \\
\Pi(\alpha) & \le & \lambda K \alpha^d \label{e:claire},
 \end{eqnarray}
where $ \epsilon $, defined by \eqref{e:application-epsilon}, tends to $0$ at $\infty$.
\end{prop}
\paragraph{Proof.} Inequalities \eqref{e:florence} and \eqref{e:claire} are part of Proposition 3.1 in \cite{G-perco-boolean-model}.
Let us prove \eqref{e:olivier}.
Let $A$ be a finite subset of $S(1)$ such that $S(1)$ is covered by the union of the balls $B(a,1/10)$, $a \in A$.
Let $K_1=K_1(d)$ denote the cardinality of $A$.

Let $\alpha > 0$. The first inequality is straightforward as 
$$
G(0,\alpha) \subset \{S(\alpha) \as S(2\alpha)\}.
$$
By definition of $A$, $S(\alpha)$ is covered by the union of the balls $B(\alpha a,\alpha/10)$, $a \in A$.
Moreover, all the balls $B(\alpha a,8\alpha/10)$, $a \in A$, are contained in $B(0,2\alpha)$.
Therefore
$$
\{S(\alpha) \as S(2\alpha) \} \subset \bigcup_{a \in A} \{S(a,\alpha/10) \as S(a,8\alpha/10) \}.
$$
By union bound, by stationarity and by definition of $K_1$ we get
$$
P(S(\alpha) \as S(2\alpha)) \le K_1 P(S(\alpha/10) \as S(8\alpha/10)).
$$
Note
\begin{equation}\label{e:telephone}
\{S(\alpha/10) \as S(8\alpha/10) \} \subset G(0,\alpha/10) \cup H(\alpha/10).
\end{equation}
Indeed, assume the existence of a path in $\Sigma$ from $S(\alpha/10)$ to $S(8 \alpha/10)$.
We can assume that the path is in $B(8\alpha/10)$.
If moreover $H(\alpha/10)$ does not occur, then  the path is in $\Sigma(B(0,10\alpha/10))$ and thus $G(0,\alpha/10)$ occurs.
Therefore,
$$
P(S(\alpha) \as S(2\alpha)) \le K_1 \Pi(\alpha/10) + K_1 P(H(\alpha/10)).
$$
But there exists a constant $K_2=K_2(d)$ such that, 
$$
P(H(\alpha/10)) \le \lambda K_2 \epsilon(\alpha/10).
$$
This is Lemma 3.4 in \cite{G-perco-boolean-model}.
The lemma follows. \finpreuve

\bigskip

Inequality \eqref{e:florence} yieds the following result.

\begin{lemma} \label{l:localisation} We use the constant $K$ from the previous lemma. Let $M>0$. Assume
\begin{equation}\label{e:controle-epsilon}
\lambda K^2 \epsilon(M) \le \frac 1 4
\end{equation}
and, for all $\alpha \in [M,10M]$,
 \begin{equation}\label{e:controle-pi}
 K\Pi(\alpha) \le \frac 1 2.
\end{equation}
Then $\lim_{\alpha \to \infty}\Pi(\alpha) = 0$.
\end{lemma}
\paragraph{Proof.} This is a consequence of \eqref{e:florence} and Lemma 3.7 in \cite{G-perco-boolean-model}. 
Showing how to apply Lemma 3.7 would not be much shorter than adapting the proof in our context. 
Therefore we choose to give a full proof.
By \eqref{e:florence}, for all $\alpha > 0$,
\begin{equation}\label{e:florence2}
K\Pi(10\alpha) \le \big( K \Pi(\alpha)\big)^2 + \lambda K^2\epsilon(\alpha).
\end{equation}
As $\epsilon$ is non-increasing, \eqref{e:controle-epsilon} yields, for all $\alpha \ge M$,
$$
K\Pi(10\alpha) \le \big( K \Pi(\alpha)\big)^2 + \frac 1 4.
$$
Therefore, if moreover $K\Pi(\alpha)\le 1/2$, then $K\Pi(10\alpha) \le 1/2$.
Using \eqref{e:controle-pi} and induction we deduce, for all $\alpha \ge M$,
$K\Pi(\alpha) \le 1/ 2$.
As a consequence,
$$
\limsup_{\alpha \to \infty} K\Pi(\alpha) \le \frac 1 2.
$$
As $\epsilon$ tends to $0$ at $\infty$ we get, using \eqref{e:florence2},
$$
\limsup_{\alpha \to \infty} K\Pi(\alpha) \le \left(\limsup_{\alpha \to \infty} K\Pi(\alpha)\right)^2.
$$
As a consequence of the two previous inequalities we get $\limsup_{\alpha \to \infty} K\Pi(\alpha)=0$. \finpreuve

\paragraph{Proof of Theorem \ref{t:explicitons}.} By \eqref{e:olivier}, 
$$
I=\{\lambda>0 : \lim_{r\rightarrow \infty}P(S(r) \as S(2r)) = 0 \} = \{\lambda > 0 : \lim_{\alpha \rightarrow \infty}\Pi(\alpha) = 0\}.
$$
Let us first prove that $I$ is non empty. 
Set $M=1$.
By \eqref{e:claire}, for small enough $\lambda > 0$, Assumptions \eqref{e:controle-epsilon} and \eqref{e:controle-pi} hold for every $\alpha \in [1,10]$.
By Lemma \ref{l:localisation}, all such $\alpha$ belong to $I$.

Let us now prove that $I$ is open. 
If $\lambda$ belongs to $I$, then any smaller positive real number belongs to $I$.
Therefore, we only have to show that, for any $\lambda \in I$, there exists $\eta > 0$ such that $\lambda+\eta \in I$.
We now fix $\lambda \in I$. 
Note that $\epsilon$ and $K$ does not depend on the density $\lambda$.
We emphasize the dependence of $\Pi$ on $\lambda$ by writing $\Pi_\lambda$.
As $\Pi_\lambda$ and $\epsilon$ tends to $0$ at infinity we can fix $M>0$ such that
\begin{equation}\label{e:seminaire_pourri}
(\lambda+1) K^2 \epsilon(M) \le \frac 1 4
\end{equation}
and, for all $\alpha \in [M,10M]$,
$$
K\Pi_\lambda(\alpha) \le \frac 1 4.
$$
Let $\eta>0$.
Consider a Boolean model $\Sigma'$ with parameters $(\eta,\nu,d)$ independent of $\Sigma$.
Then $\Sigma\cup \Sigma'$ is a Boolean model with parameters $(\lambda+\eta,\nu,d)$.
Therefore, for any $\alpha \in [M,10M]$,
\begin{eqnarray*}
\Pi_{\lambda+\eta}(\alpha) 
 & \le & \Pi_\lambda(\alpha) + P\big(\mbox{one of the random balls of }\Sigma'\mbox{ is centered in }B(0,10\alpha)\big) \\
 & \le & \Pi_\lambda(\alpha) + \eta v_d (10\alpha)^d \\
 & \le & \Pi_\lambda(\alpha) + \eta v_d (100M)^d
\end{eqnarray*}
where $v_d$ denotes the volume of the unit ball of $\R^d$.
As a consequence, for $\eta \in (0,1)$ small enough, for any $\alpha \in [M,10M]$,
$$
K\Pi_{\lambda+\eta}(\alpha) \le \frac 1 2.
$$
From \eqref{e:seminaire_pourri} and $\eta \le 1$ we also get
$$
(\lambda+\eta) K^2 \epsilon(M) \le \frac 1 4.
$$
By Lemma \ref{l:localisation}, we then deduce the convergence of $\Pi_{\lambda+\eta}(\alpha)$ to $0$ as $\alpha$ tends to $\infty$.
In other words, $\lambda+\eta$ belongs to $I$. \finpreuve

\section{Asymptotic behaviour of $T(x)/\|x\|$}

\label{s:preuve-forme}

In this section we prove Theorem \ref{t:forme}. It can be deduced from Theorem 1 in Ziesche \cite{Ziesche2}, but we give a proof of this result for self-containedness. 
The plan of proof is standard and the proof is actually particularly simple thanks to good upper bounds on $T$.
Let us first state Kingman's theorem. We choose to state this theorem as Kesten in \cite{Kesten-saint-flour} (see Theorem 2.1), following the statement proposed by Liggett in \cite{Liggett85}.

\begin{theorem} \label{th-kingman}
Suppose $(X_{m,n}, 0\le m<n)$ ($m$ and $n$ are integer) is a family of random variables satisfying:
\begin{enumerate}
\item For all integers $m,n$ such that $0 <m<n$, one has $X_{0,n}\le X_{0,m}+X_{m,n}$,
\item For each $m\geq 0$, the distribution of $(X_{m+h,m+h+k}, k\geq 1)$ does not depend on the integer $h\geq0$, 
\item For each $k\geq 1$, the sequence $(X_{nk,(n+1)k}, n\ge 0)$ is stationary and ergodic,
\item $E(X_{0,1}^+)<\infty$ and there exists a real $c$ such that, for all natural integer $n$, one has
$E(X_{0,n})\ge -cn$.
\end{enumerate}
Then
$$ \lim_{n\rightarrow \infty} \frac{X_{0,n}}{n} = \gamma   \hbox{ a.s. and in }L^1$$
where $\gamma$ is the finite constant defined by
$$\gamma = \inf_n  \frac{E (X_{0,n})}{n} .$$
\end{theorem}

Let $x \in S(1)$. 
We apply Kingman's theorem to the family defined by $X_{m,n}=T(mx,nx)$. 
\begin{itemize}
 \item For any $a,b,c \in \R^d$, $T(a,c) \le T(a,b) + T(b,c)$. This follows from the fact that the concatenation of a path from $a$ to $b$ and a path
 from $b$ to $c$ is a path from $a$ to $c$. Therefore the first assumption of Kingman's theorem holds.
 \item The process $X$ is stationary and ergodic under the action of spatial translations. Therefore the second and third assumptions of Kingman's theorem hold.
 \item For any $a,b \in \R^d$, 
 \begin{equation} \label{e:controle-T}
 0 \le T(a,b) \le \|b-a\|.
 \end{equation}
Therefore the forth assumption holds.
\end{itemize}
Thus,
$$
\lim_{n\rightarrow \infty} \frac{T(0,nx)}n = \mu(x) \mbox{ a.s.\ and in }L^1.
$$
By isotropy of the model, we get that $\mu$ does not depend on $x \in S(1)$. Therefore we drop the dependence on $x$ and write $\mu$.
We have proven
\begin{equation}\label{e:consequence-kingman}
\mbox{for all } x \in S(1), \; \lim_{n\rightarrow \infty}\frac{T(0,nx)}n = \mu(x) \mbox{ a.s. and in }L^1.
\end{equation}

\bigskip

Now we prove the uniformity of the convergence.
For any real $u$, we denote by $\lfloor u \rfloor$ its integer part and by $\{u\}$ its fractional part.
In particular, $u = \lfloor u \rfloor + \{u\}$. Let $\epsilon > 0$. Let $A$ be a finite subset of $S(1)$ such that any point of $S(1)$ is at most at distance $\epsilon$ of some point of $A$.
By \eqref{e:consequence-kingman}, with probability one, there exists $N$ such that for any $n \ge N$ and for any $x \in A$,
$$
\left| \frac{T(0,nx)}{n} - \mu \right| \le \epsilon.
$$
Let $y \in \R^d\setminus \{0\}$. Write 
\begin{equation}\label{e:ateliers-pedagogiques}
\widehat y = \frac{1}{\|y\|} y \mbox{ and } n(y)= \lfloor \|y\|\rfloor.
\end{equation}
We assume $n(y) \ge N$ and $n(y)\epsilon \ge 1$.
Let $x \in A$ be such that $\|\widehat y-x\|\le \epsilon$.
By triangle inequality for $T$ and by \eqref{e:controle-T} we get
\begin{eqnarray*}
|T(0,y)-T(0,n(y)x)| 
 & \le &  \left|T(0,y)-T\big(0, n(y) \widehat y\big)\right| + \left|T\big(0, n(y) \widehat y\big) - T(0, n(y)x\big)\right|\\
 & \le & \| y - n(y) \widehat y\| + \|n(y) \widehat y -n(y)x \| \\
  & \le & 1 + n(y) \epsilon \\
 & \le & 2n(y)\epsilon.
\end{eqnarray*}
Moreover, as $n(y) \ge N$, we also have $|T(0,n(y)x) - n(y)\mu| \le n(y)\epsilon$.
Therefore
$$
\frac{|T(0,y)-n(y)\mu|}{n(y)} \le 3\epsilon.
$$
The almost sure convergence in Theorem \ref{t:forme} follows. The convergence in $L^1$ is a straightforward consequence of the a.s. convergence and the dominated convergence with the domination $T(0,y) / \|y\| \leq 1$ for every $y \in \R^d\setminus \{0\}$.

%

\def\cprime{$'$} \def\cprime{$'$}

\end{document}